\documentclass{article}%
\usepackage{amsmath}%
\usepackage{amsfonts}%
\usepackage{amssymb}%
\usepackage{graphicx}
\newtheorem{theorem}{Theorem}

\begin{document}

\title{The ancient Greeks present: Rational Trigonometry}
\author{N J Wildberger\\School of Mathematics and Statistics \\UNSW Sydney 2052 Australia\\webpages: http://web.maths.unsw.edu/\symbol{126}norman/}
\maketitle

\begin{abstract}
Pythagoras' theorem, the area of a triangle as one half the base times the
height, and Heron's formula are amongst the most important and useful results
of ancient Greek geometry. Here we look at all three in a new and improved
light, using \textit{quadrance} not \textit{distance}. This leads to a simpler
and more elegant trigonometry, in which \textit{angle} is replaced by
\textit{spread}, and which extends to arbitrary fields and more general
quadratic forms.

\end{abstract}

\section*{Three ancient Greek theorems}

There are three classical theorems about triangles that every mathematics
student traditionally meets. To state these, consider a triangle
$\overline{A_{1}A_{2}A_{3}}$ with side lengths $d_{1}\equiv\left\vert
A_{2},A_{3}\right\vert $, $d_{2}\equiv\left\vert A_{1},A_{3}\right\vert $ and
$d_{3}\equiv\left\vert A_{1},A_{2}\right\vert $.

\begin{description}
\item[Pythagoras' theorem] \textit{The triangle }$\overline{A_{1}A_{2}A_{3}}%
$\textit{\ has a right angle at }$A_{3}$\textit{\ precisely when}%
\[
d_{1}^{2}+d_{2}^{2}=d_{3}^{3}.
\]

\item[Area of triangle] \textit{The area of a triangle is one half the length
of the base times the height.}

\item[Heron's formula] \textit{If }$s\equiv\left(  d_{1}+d_{2}+d_{3}\right)
/2$\textit{\ is the semi-perimeter of a triangle, then its area is}%
\[
\mathrm{area}=\sqrt{s\left(  s-d_{1}\right)  \left(  s-d_{2}\right)  \left(
s-d_{3}\right)  }.
\]

\end{description}

In this paper we will recast all three in simpler and more general forms. As a
reward, we find that \textit{rational trigonometry} falls into our laps,
essentially for free. Our reformulation works over a general field (not of
characteristic two), in arbitrary dimensions, and even with an arbitrary
quadratic form---see \cite{Wild1}, \cite{Wild2} and \cite{Wild3}.

\section*{Pythagoras' theorem}

Euclid and other ancient Greeks regarded \textit{area}, not \textit{distance},
as the fundamental quantity in planar geometry. Indeed they worked with a
straightedge and compass in their constructions, not a ruler and protractor. A
line segment was measured by constructing a square on it, and determining the
area of that square. Two line segments were considered equal if they were
\textit{congruent}, but this was independent of a direct notion of distance measurement.

Area is an \textit{affine concept}: more precisely \textit{proportions}
between areas are maintained by linear transformations. Even with a different
metrical geometry, for example a relativistic geometry in which $x^{2}-y^{2}$
plays the role of $x^{2}+y^{2}$, the notion of signed area defined by a
determinant applies. In other words, in planar geometry area is a basic notion
and may be considered prior to any theory of linear measurement.

To Euclid, Pythagoras' theorem is a relation about the areas of squares built
on each of the sides of a right triangle. This insight has largely been lost
in the modern formulation, but with a sheet of graph paper it is still an
attractive way to introduce students to the subject, as the area of many
simple figures can be computed by subdividing, translating and counting cells.%
%TCIMACRO{\FRAME{fhFU}{7.1413cm}{6.1734cm}{0pt}{\Qcb{Pythagoras' theorem using
%areas}}{\Qlb{Pythagoras' theorem}}{pyth4,2bw.eps}%
%width 7.1413cm;  height 6.1734cm;  depth 0pt;  original-width 7.4017cm;
%original-height 6.3917cm;  cropleft "0";  croptop "1";  cropright "1";
%cropbottom "0";  filename 'Pyth4,2BW.eps';file-properties "XNPEU";}}}%
%BeginExpansion
\begin{figure}
[h]
\begin{center}
\includegraphics[
height=6.1734cm,
width=7.1413cm
]%
{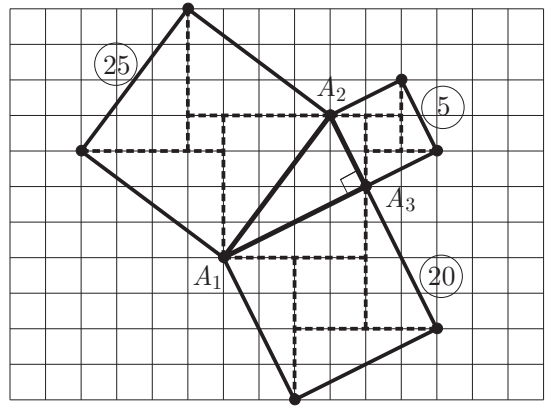}%
\caption{Pythagoras' theorem using areas}%
\label{Pythagoras' theorem}%
\end{center}
\end{figure}
%EndExpansion

The squares on the sides of triangle $\overline{A_{1}A_{2}A_{3}}$ shown in
Figure \ref{Pythagoras' theorem} have areas $5,20$ and $25.$ The largest
square for example can be seen as four triangles which can be rearranged to
get two $3\times4$ rectangles, together with a $1\times1$ square, for a total
area of $25.$ So from this point of view Pythagoras' theorem is a result which
can be established by \textit{counting}, and the use of irrational numbers to
describe lengths is not necessary. This applies to any right triangle with
rational coordinates.

Following the Greek terminology of `quadrature', we define the
\textbf{quadrance} $Q$ of a line segment to be the area of the square
constructed on it. Pythagoras' theorem allows us to assert that if
$A_{1}\equiv\left[  x_{1},y_{1}\right]  $ and $A_{2}\equiv\left[  x_{2}%
,y_{2}\right]  ,$ then the quadrance between $A_{1}$ and $A_{2}$ is%
\[
Q\left(  A_{1},A_{2}\right)  =\left(  x_{2}-x_{1}\right)  ^{2}+\left(
y_{2}-y_{1}\right)  ^{2}.
\]
So for example the quadrance between the points $\left[  0,0\right]  $ and
$\left[  1,2\right]  $ is $Q=5.$ The usual distance between the points is the
`square root' of the quadrance and requires a prior theory of irrational
numbers. Clearly the irrational number $\sqrt{5}\approx\allowbreak
2.\,\allowbreak236\,067\,977\ldots\allowbreak$ is a far more sophisticated and
complicated object than the natural number $5.$

In statistics, variance is more natural than standard deviation. In quantum
mechanics, wave functions are more basic than probability amplitudes. In
harmonic analysis, $L^{2}$ is more pleasant than $L^{1}.$ In geometry,
\textit{quadrance is more fundamental than distance}.

For a triangle $\overline{A_{1}A_{2}A_{3}}$ we define the quadrances
$Q_{1}=Q\left(  A_{2},A_{3}\right)  $, $Q_{2}=Q\left(  A_{1},A_{3}\right)  $
and $Q_{3}=Q\left(  A_{1},A_{2}\right)  $. Here then is\textit{\ Pythagoras'
theorem as the Greeks viewed it}---and it extends to arbitrary fields, to many
dimensions, and even with general quadratic forms.

\begin{theorem}
[Pythagoras]The lines $A_{1}A_{3}$ and $A_{2}A_{3}$ of the triangle
$\overline{A_{1}A_{2}A_{3}}$ are perpendicular precisely when%
\[
Q_{1}+Q_{2}=Q_{3}.
\]

\end{theorem}

\section*{Area of a triangle}

\textit{The area of a triangle is one-half the base times the height}. Let's
see how we might remove some of the irrationalities that occur with this
ancient formula by looking at an example.%
%TCIMACRO{\FRAME{fhFU}{1.9433in}{1.5338in}{0pt}{\Qcb{A triangle and an
%associated parallelogram}}{\Qlb{Euclid1}}{euclid2.eps}%
%width 1.9433in;  height 1.5338in;  depth 0pt;  original-width 1.9432in;
%original-height 1.5275in;  cropleft "0";  croptop "1";  cropright "1";
%cropbottom "0";  filename 'Euclid2.eps';file-properties "XNPEU";}} }%
%BeginExpansion
\begin{figure}
[h]
\begin{center}
\includegraphics[
height=1.5338in,
width=1.9433in
]%
{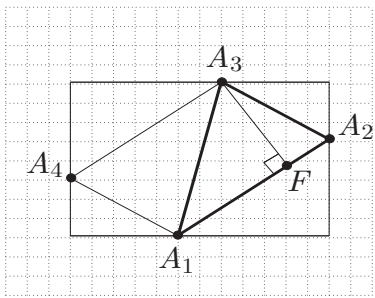}%
\caption{A triangle and an associated parallelogram}%
\label{Euclid1}%
\end{center}
\end{figure}
%EndExpansion
The area of the triangle $\overline{A_{1}A_{2}A_{3}}$ in Figure \ref{Euclid1}
is one half of the area of the associated parallelogram $\overline{A_{1}%
A_{2}A_{3}A_{4}}.$

The latter area may be calculated by removing from the circumscribed
$12\times8$ rectangle four triangles, which can be combined to form two
rectangles, one $5\times3$ and the other $7\times5.$ The area of
$\overline{A_{1}A_{2}A_{3}}$ is thus $23.$

To apply the one-half base times height rule, the base $\overline{A_{1}A_{2}}
$ by Pythagoras has length
\[
d_{3}=\left\vert A_{1},A_{2}\right\vert =\sqrt{7^{2}+5^{2}}=\sqrt{74}%
\approx8.\,\allowbreak602\,325\,267\,\allowbreak04\ldots.
\]
To find the length $h$ of the altitude $\overline{A_{3}F},$ set the origin to
be at $A_{1},$ then the line $A_{1}A_{2}$ has Cartesian equation $5x-7y=0$
while $A_{3}=\left[  2,8\right]  .$ A well-known result from coordinate
geometry then states that the distance $h=\left\vert A_{3},F\right\vert $ from
$A_{3}$ to the line $A_{1}A_{2}$ is%
\[
h=\frac{\left\vert 5\times2-7\times8\right\vert }{\sqrt{5^{2}+7^{2}}}%
=\frac{46}{\sqrt{74}}\allowbreak\approx5.\,\allowbreak
347\,391\,382\,\allowbreak22\ldots.
\]
If an engineer doing this calculation works with the surd forms of both
expressions, she will notice that the two occurrences of $\sqrt{74}$
conveniently cancel when she takes one half the product of $d_{3}$ and $h$,
giving an area of $23.$ However if she works immediately with the decimal
forms, she may be surprised that her calculator gives%
\[
\mathrm{area}\approx\frac{8.\,\allowbreak602\,325\,267\,\allowbreak
04\ldots\times~5.\,\allowbreak347\,391\,382\,\allowbreak22\ldots}{2}%
\approx23.000\,000\,000\,\allowbreak01.
\]
The usual formula forces us to descend to the level of irrational numbers and
square roots, even when the eventual answer is a natural number, and this
introduces unnecessary approximations and inaccuracies into the subject. It is
not hard to see how the use of quadrance allows us to reformulate the result.

\begin{theorem}
[Triangle area]The square of the area of a triangle is one-quarter the
quadrance of the base times the quadrance of the corresponding altitude.
\end{theorem}

As a formula, this would be%
\[
\mathrm{area}^{2}=\frac{Q\times H}{4}
\]
where $Q$ is the quadrance of the base and $H$ is the quadrance of the
altitude to that base.

\section*{Heron's or Archimedes' Theorem}

The same triangle $\overline{A_{1}A_{2}A_{3}}$ of the previous section has
side lengths
\[%
\begin{tabular}
[c]{lllll}%
$d_{1}=\sqrt{34}$ &  & $d_{2}=\sqrt{68}$ &  & $d_{3}=\sqrt{74}.$%
\end{tabular}
\]
The semi-perimeter $s,$ defined to be one half of the sum of the side lengths,
is then%
\[
s=\frac{\sqrt{34}+\sqrt{68}+\sqrt{74}}{2}\allowbreak\approx11.\,\allowbreak
339\,744\,206\,\allowbreak6\ldots.
\]

Using the usual Heron's formula, a computation with the calculator shows that
\[
\mathrm{area}=\sqrt{s\left(  s-\sqrt{34}\right)  \left(  s-\sqrt{68}\right)
\left(  s-\sqrt{74}\right)  }\approx23.000\,000.
\]
Again we have a formula involving square roots in which there appears to be a
surprising integral outcome. Let's now give another form of Heron's formula,
with a new name. Arab sources suggest that Archimedes knew Heron's formula
earlier, and the greatest mathematician of all time deserves credit for more
than he currently gets.

\begin{theorem}
[Archimedes]The area of a triangle $\overline{A_{1}A_{2}A_{3}}$ with
quadrances $Q_{1},Q_{2}$ and $Q_{3}$ is given by%
\[
16~\mathrm{area}^{2}=\left(  Q_{1}+Q_{2}+Q_{3}\right)  ^{2}-2\left(  Q_{1}%
^{2}+Q_{2}^{2}+Q_{3}^{2}\right)  .
\]

\end{theorem}

In our example the triangle has quadrances $34,68$ and $74,$ each obtained by
Pythagoras' theorem. So Archimedes' theorem states that%
\[
16~\mathrm{area}^{2}=\left(  34+68+74\right)  ^{2}-2\left(  34^{2}%
+68^{2}+74^{2}\right)  =8464
\]
and this gives an area of $23.$ In rational trigonometry, the quantity
\[
\mathcal{A}=\left(  Q_{1}+Q_{2}+Q_{3}\right)  ^{2}-2\left(  Q_{1}^{2}%
+Q_{2}^{2}+Q_{3}^{2}\right)
\]
is the \textbf{quadrea} of the triangle, and turns out to be the single most
important number associated to a triangle. Note that%
\begin{align*}
\mathcal{A} &  =4Q_{1}Q_{2}-\left(  Q_{1}+Q_{2}-Q_{3}\right)  ^{2}\\
&  =-%
\begin{vmatrix}
0 & Q_{1} & Q_{2} & 1\\
Q_{1} & 0 & Q_{3} & 1\\
Q_{2} & Q_{3} & 0 & 1\\
1 & 1 & 1 & 0
\end{vmatrix}
.
\end{align*}

It is instructive to see how to go from Heron's formula to Archimedes'
theorem. In terms of the side lengths $d_{1},d_{2}$ and $d_{3}$:%
\begin{align*}
16~\mathrm{area}^{2}  & =\left(  d_{1}+d_{2}+d_{3}\right)  \left(  d_{1}%
+d_{2}-d_{3}\right)  \left(  -d_{1}+d_{2}+d_{3}\right)  \left(  d_{1}%
-d_{2}+d_{3}\right) \\
& =\left(  \left(  d_{1}+d_{2}\right)  ^{2}-d_{3}^{2}\right)  \left(
d_{3}^{2}-\left(  d_{1}-d_{2}\right)  ^{2}\right) \\
& =\left(  2d_{1}d_{2}+\left(  d_{1}^{2}+d_{2}^{2}-d_{3}^{2}\right)  \right)
\left(  2d_{1}d_{2}-\left(  d_{1}^{2}+d_{2}^{2}-d_{3}^{2}\right)  \right) \\
& =4d_{1}^{2}d_{2}^{2}-\left(  d_{1}^{2}+d_{2}^{2}-d_{3}^{2}\right)  ^{2}\\
& =4Q_{1}Q_{2}-\left(  Q_{1}+Q_{2}-Q_{3}\right)  ^{2}\\
& =\left(  Q_{1}+Q_{2}+Q_{3}\right)  ^{2}-2\left(  Q_{1}^{2}+Q_{2}^{2}%
+Q_{3}^{2}\right)  .
\end{align*}

Archimedes' theorem implies another formula of considerable importance.

\begin{theorem}
[Triple quad formula]The three points $A_{1},A_{2}$ and $A_{3}$ are collinear
precisely when%
\[
\left(  Q_{1}+Q_{2}+Q_{3}\right)  ^{2}=2\left(  Q_{1}^{2}+Q_{2}^{2}+Q_{3}%
^{2}\right)  .
\]

\end{theorem}

The proof is of course immediate, as collinearity is equivalent to the area of
the triangle being zero.

\section*{Spread between lines}

An \textit{angle }is the ratio of a \textit{circular distance} to a
\textit{linear distance,} and this is a\textit{\ }complicated concept. To
define an angle properly \textit{you require calculus}, an important point
essentially understood by Archimedes. Vagueness about angles, and the
accompanying ambiguities in the definition of the circular functions
$\cos\theta,$ $\sin\theta$ and $\tan\theta$ weaken most calculus texts.

There is a reason that classical trigonometry is painful to students---it is
\textit{based on the wrong notions}. As a result, mathematics teachers are
forced to continually rely on $90-45-45$ and $90-60-30$ triangles for examples
and test questions, which makes the subject very narrow and repetitive.

Rational trigonometry, developed in \cite{Wild1}, see also \cite{Wild3}, shows
how to simplify and enrich the subject, leading to greater accuracy and
quicker computations. We want to show that the basic results of this new
theory \textit{follow naturally from the above presentation of Pythagoras'
theorem, the Triangle area theorem, and Archimedes' theorem}.

The key innovation is to replace angle with a completely algebraic concept.
The separation between lines $l_{1}$ and $l_{2}$ is captured rather by the
notion of \textit{spread}, which may be defined as the ratio of two quadrances
as follows.

Suppose $l_{1}$ and $l_{2}$ intersect at the point $A.$ Choose a point $B\neq
A$ on one of the lines, say $l_{1},$ and let $C$ be the foot of the
perpendicular from $B$ to $l_{2}$.%
%TCIMACRO{\FRAME{fhFU}{1.5244in}{1.3379in}{0pt}{\Qcb{Spread $s$ between two
%lines $l_{1}$ and $l_{2}$}}{\Qlb{SpreadRatio}}{spreaddef.eps}%
%width 1.5244in;  height 1.3379in;  depth 0pt;  original-width 1.5178in;
%original-height 1.3281in;  cropleft "0";  croptop "1";  cropright "1";
%cropbottom "0";  filename 'SpreadDef.eps';file-properties "XNPEU";}} }%
%BeginExpansion
\begin{figure}
[h]
\begin{center}
\includegraphics[
height=1.3379in,
width=1.5244in
]%
{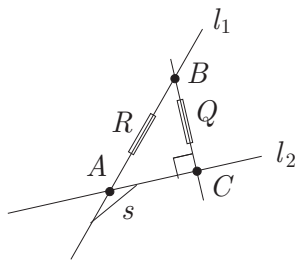}%
\caption{Spread $s$ between two lines $l_{1}$ and $l_{2}$}%
\label{SpreadRatio}%
\end{center}
\end{figure}
%EndExpansion
Then the \textbf{spread} $s$ between $l_{1}$ and $l_{2}$ is
\[
s=s\left(  l_{1},l_{2}\right)  =\frac{Q\left(  B,C\right)  }{Q\left(
A,C\right)  }=\frac{Q}{R}.
\]
This ratio is independent of the choice of $B,$ according to Thales, and is
defined \textit{between lines, not rays}. Parallel lines are defined to have
spread $s=0,$ while perpendicular lines have spread $s=1.$

The spread between the lines $l_{1}$ and $l_{2}$ with equations $a_{1}%
x+b_{1}y+c_{1}=0$ and $a_{2}x+b_{2}y+c_{2}=0$ turns out to be%
\begin{equation}
s\left(  l_{1},l_{2}\right)  =\frac{\left(  a_{1}b_{2}-a_{2}b_{1}\right)
^{2}}{\left(  a_{1}^{2}+b_{1}^{2}\right)  \left(  a_{2}^{2}+b_{2}^{2}\right)
}.\label{SpreadFormula}%
\end{equation}
Since this is a rational expression, spread becomes a useful concept also over
general fields, although in this paper we stick to the usual situation over
the decimal numbers.

You may check that the spread corresponding to $30^{\circ}$ or $150^{\circ} $
or $210^{\circ}$ or $330^{\circ}$ is $s=1/4,$ the spread corresponding to
$45^{\circ}$ or $135^{\circ}$ etc. is $s=1/2$, and the spread corresponding to
$60^{\circ}$ or $120^{\circ}$ etc. is $3/4.$

The following \textit{spread protractor} was created by M. Ossmann
\cite{Ossmann}.
%TCIMACRO{\FRAME{fhFU}{2.8755in}{1.6071in}{0pt}{\Qcb{A spread protractor}%
%}{\Qlb{SpreadProt}}{spreadprotractor.eps}%
%width 2.8755in;  height 1.6071in;  depth 0pt;  original-width 6.4446in;
%original-height 3.5908in;  cropleft "0";  croptop "1";  cropright "1";
%cropbottom "0";  filename 'SpreadProtractor.eps';file-properties "XNPEU";}} }%
%BeginExpansion
\begin{figure}
[h]
\begin{center}
\includegraphics[
height=1.6071in,
width=2.8755in
]%
{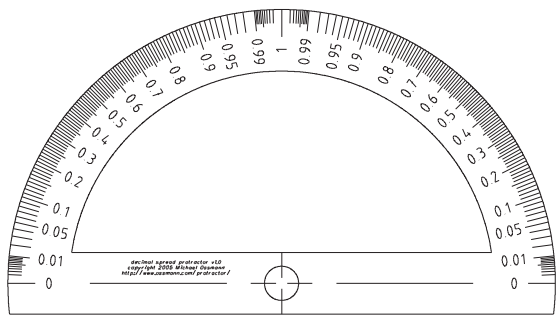}%
\caption{A spread protractor}%
\label{SpreadProt}%
\end{center}
\end{figure}
%EndExpansion

\begin{figure}
[h]
\begin{center}
\includegraphics[
height=1.1836in,
width=1.6202in
]%
{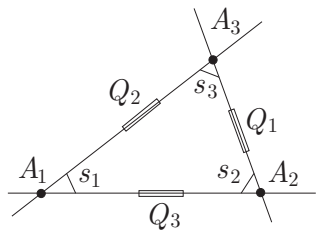}%
\caption{Quadrances and spreads of a triangle}%
\label{SpreadsQuads}%
\end{center}
\end{figure}

We use the notation that a triangle $\overline{A_{1}A_{2}A_{3}}$ has
quadrances $Q_{1},Q_{2}$ and $Q_{3},$ as well as spreads $s_{1},s_{2}$ and
$s_{3},$ labelled as in Figure \ref{SpreadsQuads}. Note the diagrammatic
conventions that help us distinguish these quantities from distance and angle.%
%TCIMACRO{\FRAME{fhFU}{1.6202in}{1.1836in}{0pt}{\Qcb{Quadrances and spreads of
%a triangle}}{\Qlb{SpreadsQuads}}{fig1.4.eps}%
%width 1.6202in;  height 1.1836in;  depth 0pt;  original-width 1.6158in;
%original-height 1.1719in;  cropleft "0";  croptop "1";  cropright "1";
%cropbottom "0";  filename 'Fig1.4.eps';file-properties "XNPEU";}} }%
%BeginExpansion
%EndExpansion

\section*{Rational trigonometry}

Let's see how to combine the three ancient Greek theorems as restated above to
derive the main laws of rational trigonometry, independent of classical
trigonometry, and without any need for transcendental functions. If $H_{3}$ is
the quadrance of the altitude from $A_{3}$ to the line $A_{1}A_{2},$ then the
Triangle area theorem and the definition of spread give%
\[
\mathrm{area}^{2}=\frac{Q_{3}\times H_{3}}{4}=\frac{Q_{3}Q_{2}s_{1}}{4}%
=\frac{Q_{3}Q_{1}s_{2}}{4}.
\]
By symmetry, we get the following analog of the Sine law.

\begin{theorem}
[Spread law]For a triangle with quadrances $Q_{1},Q_{2}$ and $Q_{3},$ and
spreads $s_{1},s_{2}$ and $s_{3},$
\[
\frac{s_{1}}{Q_{1}}=\frac{s_{2}}{Q_{2}}=\frac{s_{3}}{Q_{3}}=\frac
{4~\mathrm{area}^{2}}{Q_{1}Q_{2}Q_{3}}.
\]

\end{theorem}

By equating the formulas for $16~\mathrm{area}^{2}$ obtained from the Triangle
area theorem and Archimedes' theorem, we get%
\begin{align*}
4Q_{2}Q_{3}s_{1}  & =\left(  Q_{1}+Q_{2}+Q_{3}\right)  ^{2}-2\left(  Q_{1}%
^{2}+Q_{2}^{2}+Q_{3}^{2}\right) \\
& =2Q_{1}Q_{2}+2Q_{1}Q_{3}+2Q_{2}Q_{3}-Q_{1}^{2}-Q_{2}^{2}-Q_{3}^{2}.
\end{align*}
Rearranging gives the following analog of the Cosine law.

\begin{theorem}
[Cross law]For a triangle with quadrances $Q_{1},Q_{2}$ and $Q_{3},$ and
spreads $s_{1},s_{2}$ and $s_{3},$%
\[
\left(  Q_{1}-Q_{2}-Q_{3}\right)  ^{2}=4Q_{2}Q_{3}\left(  1-s_{1}\right)  .
\]

\end{theorem}

Now substitute $Q_{1}=s_{1}D,$ $Q_{2}=s_{2}D$ and $Q_{3}=s_{3}D,$ where
$D=Q_{1}Q_{2}Q_{3}/4~\mathrm{area}^{2}$ from the Spread law into the Cross
law, and cancel the common factor of $D^{2}$. The result is the relation%
\[
\left(  s_{1}-s_{2}-s_{3}\right)  ^{2}=4s_{2}s_{3}\left(  1-s_{1}\right)
\]
between the three spreads of a triangle, which can be rewritten more
symmetrically as follows.

\begin{theorem}
[Triple spread formula]%
\[
\left(  s_{1}+s_{2}+s_{3}\right)  ^{2}=2\left(  s_{1}^{2}+s_{2}^{2}+s_{3}%
^{2}\right)  +4s_{1}s_{2}s_{3}.
\]

\end{theorem}

This formula is a deformation of the Triple quad formula by a single cubic
term, and is the analog in rational trigonometry to the classical fact that
the three angles of a triangle sum to $\allowbreak3.\,\allowbreak
141\,592\,653\,\allowbreak59\ldots.$

The \textit{Triple quad formula}, \textit{Pythagoras' theorem}, the
\textit{Spread law}, the \textit{Cross law} and the \textit{Triple spread
formula} are the five main laws of rational trigonometry. We now see these are
closely linked to the geometrical work of the ancient Greeks.

As demonstrated at some length in \cite{Wild1}, these formulas and a few
additional ones suffice to solve the majority of trigonometric problems,
usually more simply, more accurately and more elegantly than the classical
theory involving transcendental circular functions and their inverses. As
shown in \cite{Wild2} and \cite{Wild4}, the same formulas extend to geometry
over arbitrary fields (not of characteristic two) and with general quadratic forms.

In retrospect, the blind spot first occurred with the Pythagoreans, who
initially believed that all of nature should be expressible in terms of
natural numbers and their proportions. When they discovered that the ratio of
the length of a diagonal to the length of a side of a square was the
incommensurable proportion $\sqrt{2}:1,$ legend has it that they tossed the
exposer of the secret overboard while at sea.

Had they maintained their beliefs in the workings of the Divine Mind, and
stuck with the \textit{squares of the lengths as the crucial quantities in
geometry,} then mathematics would have had a significantly different history,
Einstein's special theory of relativity would possibly have been discovered
earlier, algebraic geometry would have quite another aspect, and students
would today be studying a simpler and more elegant trigonometry---much more happily!


\begin{thebibliography}{9}                                                                                                %
\bibitem[1]{Ossmann}M. Ossmann, `Print a Protractor', download online at http://www.ossmann.com/protractor/.

\bibitem[2]{Wild1}N. J. Wildberger, \textit{Divine Proportions: Rational
Trigonometry to Universal Geometry}, Wild Egg Books, Sydney, 2005, http://wildegg.com.

\bibitem[3]{Wild2}N. J. Wildberger, Affine and Projective Rational
Trigonometry, 2006 arXiv:math/0612499.

\bibitem[4]{Wild3}N. J. Wildberger, A Rational Approach to Trigonometry,
\textit{Math Horizons}, Nov. 2007.

\bibitem[5]{Wild4}N. J. Wildberger, One dimensional metrical geometry,
\textit{Geometriae Dedicata}, \textbf{128}, no. 1, 145-166, 2007.
\end{thebibliography}
\end{document}